\documentclass{article}
\pdfoutput=1

\usepackage{graphicx}
\usepackage{amsthm}
\usepackage{amsmath, amssymb}
\usepackage{hyperref}
\usepackage[title]{appendix}

\newtheorem{theorem}{Theorem}          
\newtheorem{lemma}{Lemma}              

\begin{document}
\title{Truncated Hausdorff Moment Problem and Analytic Solution of the Time-Optimal Problem for One Class of Nonlinear Control Systems
\thanks{This work was financially supported by Polish National Science Centre grant no. 2017/25/B/ST1/01892.}}

\author{Grigory~M.~Sklyar\thanks{G.M. Sklyar is with Institute of Mathematics, University of Szczecin, Wielkopolska str. 15, Szczecin 70-451, Poland (e-mail: sklar@univ.szczecin.pl).}  \ 
		and
        Svetlana~Yu.~Ignatovich\thanks{S.Yu. Ignatovich is with V.N. Karazin Kharkiv National University, Svobody sqr. 4, Kharkiv 61022, Ukraine (e-mail: ignatovich@ukr.net).} 
}

\date{}

\maketitle

\begin{abstract}
A complete analytic solution for the time-optimal control problem for nonlinear control systems of the form $\dot x_1=u$, $\dot x_j=x_1^{j-1}$, $j=2,\ldots,n$, is obtained for arbitrary~$n$. The main goal of the paper is to present the following surprising observation: this nonlinear optimality problem leads to a truncated Hausdorff moment problem, which is applied essentially for finding the optimal time and optimal controls. 
\end{abstract}


\section{Introduction}\label{sec1}

The time-optimal control problem is one of the most well-studied topics in the optimal control theory. Among all such problems, the simplest one concerns the linear single-input system of integrators  
\begin{equation}\label{i2}
\dot x_1=u, \ \ \dot x_j=x_{j-1}, \ j=2,\ldots,n,
\end{equation}
\begin{equation}\label{i3}
x(0)=x^0, \ x(\theta)=0, \ |u(t)|\le1, \ \theta\to\min.
\end{equation}
As follows from the Pontryagin Maximum Principle 
\cite{P}, 
the optimal control takes the values $\pm1$ only and has no more than $n-1$ points of discontinuity. Therefore, in the general case, the problem \eqref{i2}, \eqref{i3} is reduced to finding $n$ numbers: the optimal time $\theta$ and switching times $t_1,\ldots,t_{n-1}$. However, except of the case $n=2$, explicit solution is not elementary. This problem was completely solved in 
\cite{KS}, \cite{KS1},
by applying ideas and technique of the Markov moment problem. More specifically, let us consider a general linear single-input time-optimal control problem
\begin{equation}\label{i0_1}
\dot x=Ax+bu,
\end{equation} 
\begin{equation}\label{i0_2}
x(0)=x^0, \ x(\theta)=0, \ |u(t)|\le1, \ \theta\to\min.
\end{equation} 
Applying the Cauchy formula for a trajectory, we obtain the following conditions for the optimal control,
\begin{equation}\label{i1}
x^0_j=\int_0^\theta g_j(t)u(t)dt, \ \ j=1,\ldots,n, \ \ |u(t)|\le1,
\end{equation} 
where $g(t)=-e^{-At}b$, which can be interpreted as a classical Markov moment problem 
\cite{M}, \cite{KN}. 
Supplemented by the optimality requirement $\theta\to\min$, the problem \eqref{i1} turns into a Markov moment \emph{min-problem} introduced in 
\cite{KS}. 
In particular, the problem \eqref{i2}--\eqref{i3} leads us to the \emph{power} Markov moment min-problem
\begin{equation}\label{i4}
\begin{array}{l}
\displaystyle (-1)^j(j-1)!x^0_{j}=\int_0^\theta t^{j-1}u(t)dt, \ \ j=1,\ldots,n,\\
 \ \ |u(t)|\le1, \ \theta\to\min.
\end{array}
\end{equation}

An idea to apply the moment problem theory to control problems was originated by N.~N.~Krasovski\u{\i}
\cite{Kr}, \cite{Kr1},
who used the results of M.~G.~Kre\u{\i}n in the abstract moment $L$-problem
\cite{AK}.
As was mentioned above, V.~I.~Korobov and G.~M.~Sklyar 
\cite{KS}, \cite{KS1}
found out that for the linear time-optimal control problem \eqref{i0_1}--\eqref{i0_2} the Markov moment problem is applicable. In 
\cite{KS} 
they obtained a complete analytic solution of the problem \eqref{i4}. The heart of the solution is as follows: \emph{For any~$n$, two special polynomials in one variable with coefficients explicitly expressed via $x^0$ are found. The maximal of roots of these polynomials coincides with the optimal time $\widehat\theta$.} After finding $\widehat\theta$, the switching times can be found as roots of another polynomial whose coefficients are explicitly expressed via  $\widehat\theta$ and $x^0$. Details of the algorithm and further references can be found in 
\cite{KSI}.

Thus, \emph{regardless of system's dimension, the solution of the time-optimal control problem \eqref{i2}--\eqref{i3} is reduced to finding roots of two polynomials in one variable.}

As was proved in 
\cite{KS1}, \cite{JMAA}, 
the integrator system \eqref{i2} approximates an arbitrary controllable linear system \eqref{i0_1} in the following sense:  after some change of variables in the system \eqref{i0_1}, optimal times and optimal controls in the problems \eqref{i2}--\eqref{i3} and \eqref{i0_1}--\eqref{i0_2} become equivalent as $x^0\to0$. Moreover, the solution of the general problem \eqref{i0_1}--\eqref{i0_2} can be found by successive solving the time-optimal control problems for the system \eqref{i2}. 

Similar analysis of the time optimality for \emph{nonlinear} systems is much more sophisticated. Let us consider the time-optimal control problem for the class of affine single-input systems, which are the most close to linear ones, 
\begin{equation}\label{i0_3}
\dot x=a(x)+b(x)u, \ \ a(0)=0,
\end{equation}
\begin{equation}\label{i0_4}
x(0)=x^0, \ x(\theta)=0, \ |u(t)|\le1, \ \theta\to\min,
\end{equation}
where  the condition $a(0)=0$ means that the origin, which is a final point, is an equilibrium of the system. We assume that the vector fields $a(x)$ and $b(x)$ are real analytic in a neighborhood of the origin. It turns out that, instead of the ``linear'' Markov moment problem \eqref{i1} with one-dimensional integrals, we obtain a ``nonlinear'' Markov moment problem 
\cite{F}, \cite{SIAM2000}
\begin{equation}\label{i0_5}
x^0=\sum v_{m_1\ldots m_k}\xi_{m_1\ldots m_k}(\theta,u), \ \ |u(t)|\le1,
\end{equation}
where $\xi_{m_1\ldots m_k}(\theta,u)$ are nonlinear functionals of $u$ of the form
\begin{equation}\label{i0_6}
\xi_{m_1\ldots m_k}(\theta,u)=\int_0^\theta\int_0^{\tau_1}\cdots\int_0^{\tau_{k-1}}\prod_{i=1}^k\tau_i^{m_i}u(\tau_i)d\tau_k\cdots d\tau_1.
\end{equation}
We refer to \eqref{i0_6} as \emph{nonlinear power moments}. Vector coefficients $v_{m_1\ldots m_k}\in{\mathbb{R}}^n$ are found by use of values of $a$ and $b$ and all their derivatives at the origin. For the linear system \eqref{i0_1}, the representation \eqref{i0_5} reads
$$x^0=\sum_{m=0}^\infty v_m\xi_m(\theta,u)=\sum_{m=0}^\infty \frac{(-1)^{m+1}}{m!}A^mb\int_0^\theta\tau^m u(\tau)d\tau,
$$
what coincides with \eqref{i1}. 

In
\cite{SIAM2000}
we described a class of nonlinear systems \eqref{i0_3} that can be approximated by linear ones in the sense of time optimality mentioned above: after some change of variables in the system \eqref{i0_3}, optimal times and optimal controls in the problems \eqref{i2}--\eqref{i3} and \eqref{i0_3}--\eqref{i0_4} become equivalent as $x^0\to0$. We called such systems \emph{essentially linear}. In a general case, as a class of simple systems that approximate arbitrary affine systems in the sense of time optimality, \emph{homogeneous systems} should be considered. For homogeneous systems, the representation \eqref{i0_5} has the form 
\begin{equation}\label{i0_7}
x^0_j=\sum_{m_1+\cdots+m_k+k=w_j} v^j_{m_1\ldots m_k}\xi_{m_1\ldots m_k}(\theta,u), \ j=1,\ldots,n,
\end{equation}
where $w_1\le\cdots\le w_n$ are some integers. One limiting case of such a representation is when the right hand side contains the unique one-dimensional integral, as in \eqref{i4}, what corresponds to the linear integrator system \eqref{i2}. The opposite limiting case is when the right hand side contains the unique multiple integral of maximal multiplicity; one can show that this is possible if
\begin{equation}\label{i0_8}
\begin{array}{l}
\displaystyle x^0_1=-\int_0^\theta u(\tau)d\tau,\\
\displaystyle x^0_{j}=(-1)^{j}(j-1)!\int_0^\theta\cdots\int_0^{\tau_{j-2}}\!\!\!\tau_{j-1}\prod_{i=1}^{j-1}u(\tau_i)d\tau_{j-1}\cdots d\tau_1,\\
 \ \ \ j=2,\ldots,n,
\end{array}
\end{equation}
what corresponds to the system
\begin{equation}\label{i0_9}
\dot x_1=u, \ \ \dot x_j=x_1^{j-1}, \ j=2,\ldots,n.
\end{equation}
The representation \eqref{i0_8} means that  \eqref{i0_9} can be considered as a ``dual-to-integrator'' system. As a homogeneous system, it approximates a certain class of affine systems in the sense of time optimality
\cite{SIAM2003}. 

In
\cite{SIS} 
we considered more general systems of the form 
\begin{equation}\label{s0}
\dot x_1=u, \ \ \dot x_j=P_j(x_1), \ j=2,\ldots,n,
\end{equation}
where $P_2(z),\ldots,P_n(z)$ are real analytic functions. We proved that the time-optimal control can take values $\pm1$ and $0$ only and has a finite number of switchings; moreover, we described a character of all time-optimal controls. 

In the present paper we consider the time-optimal control problem for the dual-to-integrator system \eqref{i0_9}. Our main goal is to present the following surprising observation: this optimality problem also leads to a moment problem, however, in this case we deal with the classical truncated \emph{Hausdorff moment problem} 
\cite{KN}, \cite{H}.
This allows us to use profound ideas and  methods of the classical moment theory and, as a result, leads to finding an analytic solution of the time-optimal control problem for the system  \eqref{i0_9}. 

The solution we have found can be summarized as follows. For any $n$, in dependence on its parity, only four or five types of optimal control are possible. For each type, in order to find the optimal time one needs to solve a system of at most two special polynomial equations in two variables with coefficients explicitly expressed via $x^0$ (in some cases, only one polynomial equation in one variable should be solved). Thus, \emph{regardless of system's dimension, the solution of the time-optimal control problem for the system \eqref{i0_9} is reduced to solving a certain polynomial equation or a certain system of two polynomial equations.}

The rest of the paper is organized as follows. In Section~\ref{sec_2} a preliminary discussion is given. In particular, we recall the results obtained in 
\cite{SIS}
for general systems of the form \eqref{s0}. In Section~\ref{subsec_3} we apply these results to the system \eqref{i0_9} and demonstrate how the truncated Hausdorff moment problem arises. Further we specify the statement of the Hausdorff moment problem for generic points and describe solvability conditions in Lemma~\ref{Lem}.  A solvability theorem for a truncated Hausdorff moment problem is recalled in Appendix~\ref{app1}; the proof of Lemma~\ref{Lem} is given in Appendix~\ref{app2}. In Section~\ref{sec_3} we formulate the main result of the paper (Theorem~\ref{Thm}) and propose an algorithm for finding optimal controls. Section~\ref{sec_4} contains several illustrative examples.
      
\section{Preliminary discussion}\label{sec_2}

We consider the following time-optimal control problem 
\begin{equation}\label{s1}
\dot x_1=u, \ \ \dot x_j=x_1^{j-1}, \ j=2,\ldots,n,
\end{equation}
\begin{equation}\label{s2}
x(0)=x^0, \ x(\theta)=0, \ |u(t)|\le1, \ \theta\to\min.
\end{equation}
The case $n=2$ is linear; it is well known 
\cite{P} 
that optimal controls take values $\pm1$ only and have no more than one switching point. The case  $n=3$ was completely solved in 
\cite{SIAM2003}. 
It was shown that optimal controls take values $\pm1$ and $0$ only. More specifically, the controllability domain is broken into two subsets with nonempty interior: the first one contains points $x^0$ which correspond to optimal controls taking values $\pm 1$ (bang-bang) and the second one contains points $x^0$ which correspond to optimal controls taking values $\pm1$ and $0$ (singular). It was shown that, outside a set of zero measure, four possible types of optimal controls are possible and domains corresponding to these types of controls are disjoint. In this paper we consider the general case $n\ge4$.

Now we recall a description of optimal controls obtained in 
\cite{SIS} 
for general systems of the form \eqref{s0}. For definiteness, suppose $x_1^0\ge0$ (the case $x_1^0\le0$ is treated completely analogously). Suppose an optimal control exists; denote it by $\widehat u(t)$. Let $\widehat x(t)$ be the corresponding optimal trajectory and $\widehat\theta$ be the optimal time.  Then, as was shown in 
\cite{SIS}, 
the optimal control $\widehat u(t)$ takes the values $\pm1$ and $0$ only and has a finite number of switching points. Moreover, there exists a polynomial in one variable $P(z)=\sum_{j=0}^{n-1}p_jz^j$ such that:

--- if $\tau$ is a switching point of $\widehat u(t)$, then $z=\widehat x_1(\tau)$ is a root of $P(z)$;

--- if $\tau$ is a switching point such that $\widehat u(\tau-0)=0$ or $\widehat u(\tau+0)=0$, then  $z=\widehat x_1(\tau)$ is a multiple root of $P(z)$;

--- $P(\widehat x_1(t))\ge0$ for $t\in[0,\widehat\theta]$, i.e., $\widehat x_1(t)$ belongs to that connected component of the set $\{z:P(z)\ge0\}$ which contains the point $z=0$.

Moreover, as was shown in  
\cite{SIS}, 
an optimal control $\widehat u(t)$ can be chosen in the ``stair-step form'' 
\begin{equation}\label{s3}
\widehat u(t)=(p_{a_1}\circ n_{\sigma_1}\circ m_{b_1}\circ\cdots\circ n_{\sigma_k}\circ m_{b_k}\circ n_{\sigma_{k+1}}\circ p_{a_2})(t),
\end{equation}
where  the following notation is used
$$p_a(t)\equiv 1, \ \ m_a(t)\equiv -1, \ \ n_a(t)\equiv 0, \ \ t\in[0,a],
$$
$k\ge1$, $a_1,a_2\ge0$, $b_1,\ldots,b_k>0$, $\sigma_1\ge0$, $\sigma_{k+1}\ge0$, $\sigma_2,\ldots,\sigma_{k}>0$, and $\circ$ means a concatenation, i.e., for functions $\varphi^1(t)$, $t\in[0,t_1]$, and $\varphi^2(t)$, $t\in[0,t_2]$, one has
$$
(\varphi^1\circ \varphi^2)(t)=\left\{\begin{array}{l}
  \varphi^1(t)\ \mbox{ for   } t\in[0,t_1), \\
  \varphi^2(t-t_1)\ \mbox{ for }\ t\in[t_1,t_1+t_2].
\end{array}\right.
$$
Additionally,  
\begin{equation}\label{s3_1}
\mbox{if } \ a_2=0, \ \mbox{ then } \ \sigma_{k+1}=0,
\end{equation}
$x_1^0+a_1+a_2=b_1+\cdots+b_k$, and 
\begin{equation}\label{s3_01}
x_1^0+a_1\not=b_1+\cdots+b_j, \ \ j=1,\ldots,k-1.
\end{equation}

\section{Time-optimal control problem as a Hausdorff moment problem}\label{subsec_3}

Now we pass to the problem \eqref{s1}--\eqref{s2}. Substituting the optimal control \eqref{s3} to the Cauchy problem $\dot x_1=u(t)$, $x_1(0)=x_1^0$, we obtain the first component of the optimal trajectory $\widehat x_1(t)$. We get
\begin{equation}\label{s4}
\widehat x_1(t)=(i^{z_0}_{a_1}\circ c^{z_1}_{\sigma_1}\circ d^{z_1}_{b_1}\circ\cdots\circ c^{z_k}_{\sigma_k}\circ d^{z_k}_{b_k}\circ c^{z_{k+1}}_{\sigma_{k+1}}\circ i^{z_{k+1}}_{a_2})(t),
\end{equation}
where 
$$i^z_a(t)=t+z, \ \ d^z_a(t)=-t+z, \ \ c^z_a(t)\equiv z, \ \ t\in[0,a],
$$
and $z_0=x_1^0$, $z_1=z_0+a_1$, $z_{j+1}=z_j-b_j$, $j=1,\ldots,k$. In particular, \eqref{s3_01} implies that $z_2,\ldots,z_{k}\not=0$. Besides,
\begin{equation}\label{s4_1}
\widehat \theta=a_1+a_2+\sum_{s=1}^kb_s+\sum_{s=1}^{k+1}\sigma_s.
\end{equation}

Now let us turn to the rest differential equations in \eqref{s1}. Substituting $\widehat x_1(t)$ to the Cauchy problems $\dot x_j=\widehat x_1^{j-1}(t)$, $x_j(0)=x_j^0$ and taking into account end conditions $x_j(\widehat\theta)=0$, we get
\begin{equation}\label{s4_0}
-x_{j}^0=\int_0^{\widehat \theta}\widehat x_1^{j-1}(t)dt, \ \ j=2,\ldots,n,
\end{equation} 
where $\widehat x_1(t)$ has the form \eqref{s4}.
Since  for functions $\varphi^1(t)$, $t\in[0,t_1]$, and $\varphi^2(t)$, $t\in[0,t_2]$, we have
$$\int_0^{t_1+t_2}(\varphi^1\circ \varphi^2)(t)dt=\int_0^{t_1}\varphi^1(t)dt+\int_0^{t_2}\varphi^2(t)dt,$$ 
the equatilies \eqref{s4_0} can be written as 
\begin{equation}\label{s5}
\begin{array}{c}
\displaystyle
-x_{j}^0=\int_0^{a_1}(t+z_0)^{j-1}dt+\int_0^{a_2}(t+z_{k+1})^{j-1}dt+\\
\displaystyle+\sum_{s=1}^{k}\int_0^{b_s}(-t+z_{s})^{j-1}dt+\sum_{s=1}^{k+1}\int_0^{\sigma_s}z_s^{j-1}dt.\end{array}
\end{equation}

Let us simplify these expressions. First, taking into account equalities $z_0=x_1^0$ and $z_{k+1}=z_1-\sum_{s=1}^kb_s=-a_2$, we denote
\begin{equation}\label{s6}
a=z_{k+1}=-a_2\le0, \quad b=z_1=x_1^0+a_1\ge x_1^0 
\end{equation} 
and get 
$$\int_0^{a_1}(t+z_0)^{j-1}dt+\int_0^{a_2}(t+z_{k+1})^{j-1}dt=$$
$$=\int_{x_1^0}^{b}t^{j-1}dt+\int_{a}^0t^{j-1}dt=\frac{b^{j}-a^{j}-(x_1^0)^{j}}{j}.
$$
Since $z_{s+1}=z_s-b_s$, we get
$$\sum_{s=1}^{k}\int_0^{b_s}(-t+z_{s})^{j-1}dt=-\sum_{s=1}^{k}\int_{z_s}^{z_{s+1}}t^{j-1}dt=$$
$$=-\int_{z_1}^{z_{k+1}}t^{j-1}dt=\int_a^bt^{j-1}dt=\frac{b^{j}-a^{j}}{j}.
$$
Finally, observe that
$$\int_0^{\sigma_s}z_s^{j-1}dt=z_s^{j-1}\sigma_s, \ s=1,\ldots, k+1.
$$
Therefore, by denoting
\begin{equation}\label{s7}
c_{j}=-x_{j}^0+\frac{(x_1^0)^{j}-2b^{j}+2a^{j}}{j},  \ 
j=2,\ldots,n,
\end{equation} 
we transform \eqref{s5} to the following form
\begin{equation}\label{s8}
c_{j}=b^{j-1}\sigma_1+\sum_{s=2}^{k}z_s^{j-1}\sigma_s+a^{j-1}\sigma_{k+1},  \ j=2,\ldots,n.
\end{equation} 
Besides, \eqref{s4_1} gives 
\begin{equation}\label{s9_0}
\widehat \theta=2b-2a-x_1^0+\sum_{s=1}^{k+1}\sigma_s.
\end{equation} 
Then, denoting 
\begin{equation}\label{s9}
c_1=\widehat \theta+x_1^0-2b+2a,
\end{equation} 
we obtain 
\begin{equation}\label{s10}
c_1=\sum_{s=1}^{k+1}\sigma_s.
\end{equation} 

Now we explain the sense of equalities \eqref{s8}, \eqref{s10}. Imagine that $a$, $b$, and $\widehat\theta$ are already known, then \eqref{s7}, \eqref{s9} define the sequence of numbers $\{c_1,\ldots,c_{n}\}$. Let us consider the truncated Hausdorff power moment problem for the sequence $\{c_1,\ldots,c_{n}\}$ on the segment $[a,b]$: \emph{to find a non-decreasing function $\sigma(z)$ such that the following moment equalities hold} 
\begin{equation}\label{s10_1}
c_j=\int_a^bz^{j-1}d\sigma(z), \ \ j=1,\ldots, n.
\end{equation}
\emph{Equalities \eqref{s8} and \eqref{s10} mean that the moment problem \eqref{s10_1} has a solution; this solution is a step function with jumps at the points $a,z_2,\ldots,z_k,b$ and jump values $\sigma_1,\ldots,\sigma_{k+1}$ respectively.} Some key results on the truncated Hausdorff power moment problem  used below are presented in Appendix~\ref{app1}.

Let us express $x_1^0$ from \eqref{s9_0} as
$$x_1^0=-\widehat \theta+2b-2a+\sum_{s=1}^{k+1}\sigma_s
$$
and then substitute it to the right hand side of the expression of $x_{j}^0$ from \eqref{s7}, \eqref{s8}
$$
x_{j}^0=\frac{(x_1^0)^{j}-2b^{j}+2a^{j}}{j}-b^{j-1}\sigma_1-\sum_{s=2}^{k}z_s^{j-1}\sigma_s-a^{j-1}\sigma_{k+1},
$$
$j=2,\ldots,n$. Thus, $x^0$ is a polynomial vector-function of the parameters $\widehat \theta$, $a$ (if $a\not=0$), $b$ (if $b\not=x_1^0$), $\sigma_1$ and $\sigma_{k+1}$ (if they are non-zero), $\sigma_2,\ldots,\sigma_k$, and $z_2,\ldots,z_k$. If the number of parameters is less than $n$, the image of such a vector-function has zero measure. In the rest of the paper we are interested only in those points $x^0$ which are defined by no less than $n$ parameters and, additionally, are such that $x_1^0\not=0$. We refer to these points as \emph{generic}. 

It turns out that for generic points, along with condition \eqref{s3_1}, which can be rewritten as 
\begin{equation}\label{s3_3}
\mbox{if } \ a=0, \ \mbox{ then } \ \sigma_{k+1}=0,
\end{equation}
the following property hold:
\begin{equation}\label{s3_2}
\mbox{if } \ b=x^0_1, \ \mbox{ then } \ \sigma_{1}=0.
\end{equation}
If fact, suppose the contrary, i.e., let $b=x^0_1$ and $\sigma_1>0$. Due to \eqref{s3_3}, three cases are possible:

--- if $a<0$ and $\sigma_{k+1}>0$, then the polynomial $P(z)$ has at least $k+1$ multiple roots. Hence, $n-1\ge 2(k+1)$, what implies $n\ge 2k+3$. In this case we have  $2k+2$ parameters (namely, $\widehat \theta$, $a$, $z_2,\ldots,z_k$, and $\sigma_1,\ldots,\sigma_{k+1}$);  

--- if  $a<0$ and $\sigma_{k+1}=0$, then the polynomial $P(z)$ has at least $k$ multiple roots and one more (maybe single) root. Hence, $n-1\ge 2k+1$, what implies $n\ge 2k+2$. In this case we have  $2k+1$ parameters (namely, $\widehat \theta$, $a$, $z_2,\ldots,z_k$, and $\sigma_1,\ldots,\sigma_{k}$); 

--- if  $a=0$ and $\sigma_{k+1}=0$, then the polynomial $P(z)$ has at least $k$ multiple roots. Hence, $n-1\ge 2k$, what implies $n\ge 2k+1$. In this case we have  $2k$ parameters (namely, $\widehat \theta$, $a$, $z_2,\ldots,z_k$, and $\sigma_2,\ldots,\sigma_{k}$).

Thus, in any case the  number of parameters is less than $n$, what is impossible for a generic point. This proves \eqref{s3_2}.

Now we analyze the possible values of $k$ for generic points. In view of conditions \eqref{s3_3} and \eqref{s3_2}, we get three cases for $b$ and $\sigma_1$: 

--- $b=x^0_1$, $\sigma_1=0$,

--- $b>x^0_1$, $\sigma_1=0$,

--- $b>x^0_1$, $\sigma_1>0$.

\noindent Independently, three cases for $a$ and $\sigma_{k+1}$ are possible: 

--- $a=0$, $\sigma_{k+1}=0$,

--- $a<0$, $\sigma_{k+1}=0$,

--- $a<0$, $\sigma_{k+1}>0$.

Generally, we have nine cases. As an example, we analyze one of them, the rest can be considered in a similar way. Suppose $b=x^0_1$, $\sigma_1=0$ and $a=0$, $\sigma_{k+1}=0$. Then the number of parameters equals $2k-1$ (they are $\widehat \theta$, $\sigma_2,\ldots,\sigma_k$, and $z_2,\ldots,z_k$), therefore, for generic points $n\le 2k-1$. On the other hand, $P(z)$ has no less than $k-1$ multiple roots, therefore, $n-1\ge 2(k-1)$, what gives $n\ge 2k-1$. Hence, $n=2k-1$. If $n$ is written as $n=2m+1$, then $k=m+1$.

Such an analysis shows that for a generic point four cases for even $n$ and five cases for odd  $n$ are possible.  We list all these cases in Table~\ref{tab1}. Besides of information on $n$ and $k$, each cell contains the type marker and values of parameter $d$ mentioned in Lemma~\ref{Lem} below. 

For illustration, Table~\ref{tab2} contains sketches of $\widehat x_1(t)$ for all the cases when $n=4$ and $n=5$.

\begin{table}\caption{Possible Cases for Generic Points}\label{tab1}
$$
\begin{array}{c||c|c|c|}
& b=x_1^0,& b>x_1^0, & b>x_1^0, \\
& \sigma_{1}=0 & \sigma_{1}=0 &\sigma_{1}>0 \\\hline\hline\rule{0pt}{10pt}
a=0,   &\underline{\mbox{Case 1}} &\underline{\mbox{Case 2}} &\underline{\mbox{Case 3}}\\
\sigma_{k+1}=0 & n=2m+1& n=2m,& n=2m+1   \\
& k=m+1, &k=m & k=m \\
& \mbox{type }A& \mbox{type }A &\mbox{type }C \\
& d=0 & d=0,1 &d=0,1 \\[2pt]\hline\rule{0pt}{10pt}
a<0,  &\underline{\mbox{Case 4}} &\underline{\mbox{Case 5}} &\underline{\mbox{Case 6}}\\
\sigma_{k+1}=0 & n=2m&n=2m+1 &n=2m   \\
& k=m &k=m & k=m-1\\
&\mbox{type }A &\mbox{type }A &\mbox{type }C\\
&d=0,1 &d=0,1,2 & d=0,1,2\\[2pt]\hline\rule{0pt}{10pt}
a<0,  &\underline{\mbox{Case 7}} &\underline{\mbox{Case 8}} &\underline{\mbox{Case 9}}\\
\sigma_{k+1}>0 & n=2m+1 &n=2m &n=2m+1   \\
&k=m &k=m-1 & k=m-1 \\
&\mbox{type }D  &\mbox{type }D &\mbox{type }B\\
&d=0,1  &d=0,1,2 &d=0,1,2\\\hline
\end{array}
$$ 
\end{table}
\begin{table}\caption{Sketches of $\widehat x_1(t)$}\label{tab2}
$$\begin{array}{c||c|c|c|}
& b=x_1^0,& b>x_1^0, & b>x_1^0, \\
& \sigma_{1}=0 & \sigma_{1}=0 &\sigma_{1}>0 \\\hline\hline\rule{0pt}{10pt}
a=0,   &n=5 &n=4 &n=5\\
\sigma_{k+1}=0 &  
\begin{picture}(45,0)(0,20)
\qbezier(5,3)(5,20)(5,25)\qbezier(0,5)(20,5)(40,5)
\qbezier(5,20)(6,19)(10,15)\qbezier(10,15)(11,15)(17,15)
\qbezier(17,15)(18,14)(22,10)\qbezier(22,10)(26,10)(30,10)
\qbezier(30,10)(31,9)(35,5)
\end{picture} & 
\begin{picture}(45,0)(0,20)
\qbezier(5,3)(5,20)(5,25)\qbezier(0,5)(20,5)(40,5)
\qbezier(5,15)(6,16)(10,20)
\qbezier(10,20)(11,19)(20,10)\qbezier(20,10)(26,10)(30,10)
\qbezier(30,10)(31,9)(35,5)
\end{picture} & 
\begin{picture}(45,0)(0,20)
\qbezier(5,3)(5,20)(5,25)\qbezier(0,5)(20,5)(40,5)
\qbezier(5,15)(6,16)(10,20)\qbezier(10,20)(11,20)(15,20)
\qbezier(15,20)(16,19)(25,10)\qbezier(25,10)(26,10)(30,10)
\qbezier(30,10)(31,9)(35,5)
\end{picture} \\ & & & \\ & & &   \\\hline\rule{0pt}{10pt}
a<0,  &n=4 &n=5 &n=4\\
\sigma_{k+1}=0 &  
\begin{picture}(45,0)(0,20)
\qbezier(5,3)(5,20)(5,25)\qbezier(0,10)(20,10)(40,10)
\qbezier(5,20)(6,19)(10,15)\qbezier(10,15)(11,15)(20,15)
\qbezier(20,15)(25,10)(30,5)\qbezier(30,5)(31,6)(35,10)
\end{picture} & 
\begin{picture}(45,0)(0,20)
\qbezier(5,3)(5,20)(5,25)\qbezier(0,10)(20,10)(40,10)
\qbezier(5,15)(6,16)(12,22)
\qbezier(12,22)(13,21)(17,17)\qbezier(17,17)(18,17)(22,17)
\qbezier(22,17)(23,16)(34,5)\qbezier(34,5)(35,6)(39,10)
\end{picture} & 
\begin{picture}(45,0)(0,20)
\qbezier(5,3)(5,20)(5,25)\qbezier(0,10)(20,10)(40,10)
\qbezier(5,15)(6,16)(10,20)\qbezier(10,20)(11,20)(15,20)
\qbezier(15,20)(25,10)(30,5)\qbezier(30,5)(31,6)(35,10)
\end{picture} \\ & & & \\ & & &
\\\hline\rule{0pt}{10pt}
a<0,  &n=5 &n=4 &n=5\\
\sigma_{k+1}>0 &  
\begin{picture}(45,0)(0,20)
\qbezier(5,3)(5,20)(5,25)\qbezier(0,10)(20,10)(40,10)
\qbezier(5,20)(6,19)(10,15)\qbezier(10,15)(11,15)(15,15)
\qbezier(15,15)(16,14)(25,5)\qbezier(25,5)(26,5)(30,5)
\qbezier(30,5)(31,6)(35,10)
\end{picture} & 
\begin{picture}(45,0)(0,20)
\qbezier(5,3)(5,20)(5,25)\qbezier(0,10)(20,10)(40,10)
\qbezier(5,15)(6,16)(10,20)
\qbezier(10,20)(11,19)(25,5)\qbezier(25,5)(26,5)(30,5)
\qbezier(30,5)(31,6)(35,10)
\end{picture} & 
\begin{picture}(45,0)(0,20)
\qbezier(5,3)(5,20)(5,25)\qbezier(0,10)(20,10)(40,10)
\qbezier(5,15)(6,16)(10,20)\qbezier(10,20)(11,20)(15,20)
\qbezier(15,20)(16,19)(30,5)\qbezier(30,5)(31,5)(35,5)
\qbezier(35,5)(36,6)(40,10)
\end{picture} \\ & & & \\ & & &
\\\hline
\end{array}
$$ 
\end{table}


The following lemma describes necessary and sufficient conditions for solvability of a Hausdorff moment problem; below we apply them for the problem \eqref{s8}, \eqref{s10}.

\begin{lemma}\label{Lem} Suppose $n\ge4$. Let a segment $[a,b]$ be fixed and numbers $c_1,\ldots,c_n$ be given. For brevity, denote
$$c^{a}_j=c_{j+1}-ac_{j}, \ \ j=1,\ldots,n-1,$$
$$c^{b}_j=-c_{j+1}+bc_{j}, \ \ j=1,\ldots,n-1,$$ 
$$c^{a,b}_j=-c_{j+2}+(a+b)c_{j+1}-abc_{j}, \ \ j=1,\ldots,n-2.$$

$(A)$ There exists a non-decreasing function $\sigma(z)$ having exactly $k-1$ nonzero points of growing between $a$ and $b$, where $2k-1\le n$, i.e., the following  representation hold
\begin{equation}\label{L1}
c_j=\sum_{s=2}^k z_s^{j-1}\sigma_s, \ \ j=1,\ldots,n,
\end{equation}
where 
\begin{equation}\label{L3} 
b>z_2>\cdots>z_k>a,
\end{equation}
\begin{equation}\label{L4} 
z_2,\ldots,z_k\ne0,
\end{equation}
\begin{equation}\label{L2} 
\sigma_2,\ldots,\sigma_k>0,
\end{equation}
if and only if
\begin{enumerate}
\item[$(A_1)$] the matrices $\{c_{i+j-1+d}\}_{i,j=1}^k$ are singular for $0\le d\le n+1-2k$, i.e., 
\begin{equation}\label{L5} 
\det\{c_{i+j-1+d}\}_{i,j=1}^k=0, \ \ 0\le d\le n+1-2k;
\end{equation}

\item[$(A_2)$] the matrices $\{c_{i+j-1}\}_{i,j=1}^{k-1}$ and $\{c_{i+j+1}\}_{i,j=1}^{k-1}$ are positive definite;

\item[$(A_3)$] the matrix $\{c^{a,b}_{i+j-1}\}_{i,j=1}^{k-1}$ is positive definite.
\end{enumerate}

The points $z_2,\ldots,z_k$ can be found as the roots of the equation
\begin{equation}\label{L5_1}
\det\left(\begin{array}{cccc}
c_1&c_2&\cdots&c_{k}\\
\cdots&\cdots&\cdots&\cdots\\
c_{k-1}&c_{k}&\cdots&c_{2k-2}\\
1&z&\cdots&z^{k-1}
\end{array}
\right)=0.
\end{equation}
If $z_2,\ldots,z_k$ are known, the numbers $\sigma_2,\ldots,\sigma_{k}$ can be found from the first $k-1$ equations of the system \eqref{L1}.

$(B)$ There exists a non-decreasing function $\sigma(z)$ having exactly $k+1$ nonzero points of growing, including $a$ and $b$, where $2k+1\le n$, i.e., the following  representation hold
\begin{equation}\label{L6}
c_j=b^{j-1}\sigma_1+\sum_{s=2}^k z_s^{j-1}\sigma_s+a^{j-1}\sigma_{k+1}, \ \ j=1,\ldots,n,
\end{equation}
where  \eqref{L3}, \eqref{L4} and 
\begin{equation}\label{L7} 
\sigma_1,\ldots,\sigma_{k+1}>0
\end{equation}
hold if and only if
\begin{enumerate}
\item[$(B_1)$] the matrices $\{c^{a,b}_{i+j-1+d}\}_{i,j=1}^{k}$ are singular for $0\le d\le n-1-2k$, i.e., 
\begin{equation}\label{L8} 
\det\{c^{a,b}_{i+j-1+d}\}_{i,j=1}^{k}=0, \ \ 0\le d\le n-1-2k;
\end{equation}

\item[$(B_2)$] the matrices $\{c^{a,b}_{i+j-1}\}_{i,j=1}^{k-1}$ and $\{c^{a,b}_{i+j+1}\}_{i,j=1}^{k-1}$ are positive definite;

\item[$(B_3)$] the matrix $\{c_{i+j-1}\}_{i,j=1}^{k+1}$ is positive definite.
\end{enumerate}

The points $z_2,\ldots,z_k$ can be found as the roots of the equation
\begin{equation}\label{L5_2}
\det\left(\begin{array}{cccc}
c^{a,b}_1&c^{a,b}_2&\cdots&c^{a,b}_{k}\\
\cdots&\cdots&\cdots&\cdots\\
c^{a,b}_{k-1}&c^{a,b}_{k}&\cdots&c^{a,b}_{2k-2}\\
1&z&\cdots&z^{k-1}
\end{array}
\right)=0.
\end{equation}
If $z_2,\ldots,z_k$ are known, the numbers $\sigma_1,\ldots,\sigma_{k+1}$ can be found from the first $k+1$ equations of the system \eqref{L6}.

$(C)$ There exists a non-decreasing function $\sigma(z)$ having exactly $k$ nonzero points of growing, including $b$ but not including $a$, where $2k\le n$, i.e., the following  representation hold
\begin{equation}\label{L9}
c_j=b^{j-1}\sigma_1+\sum_{s=2}^k z_s^{j-1}\sigma_s, \ \ j=1,\ldots,n,
\end{equation}
where  \eqref{L3}, \eqref{L4} and 
\begin{equation}\label{L10} 
\sigma_1,\ldots,\sigma_{k}>0
\end{equation}
hold if and only if
\begin{enumerate}
\item[$(C_1)$] the matrices $\{c^{b}_{i+j-1+d}\}_{i,j=1}^{k}$ are singular for $0\le d\le n-2k$, i.e., 
\begin{equation}\label{L11} 
\det\{c^{b}_{i+j-1+d}\}_{i,j=1}^{k}=0, \ \ 0\le d\le n-2k;
\end{equation}

\item[$(C_2)$] the matrices $\{c^{b}_{i+j-1}\}_{i,j=1}^{k-1}$ and $\{c^{b}_{i+j+1}\}_{i,j=1}^{k-1}$ are positive definite;

\item[$(C_3)$] the matrix $\{c^a_{i+j-1}\}_{i,j=1}^{k}$ is positive definite.
\end{enumerate}

The points $z_2,\ldots,z_k$ can be found as the roots of the equation
\begin{equation}\label{L5_3}
\det\left(\begin{array}{cccc}
c_1^b&c_2^b&\cdots&c_{k}^b\\
\cdots&\cdots&\cdots&\cdots\\
c^{b}_{k-1}&c^{b}_{k}&\cdots&c^{b}_{2k-2}\\
1&z&\cdots&z^{k-1}
\end{array}
\right)=0.
\end{equation}
If $z_2,\ldots,z_k$ are known, the numbers $\sigma_1,\ldots,\sigma_{k}$ can be found from the first $k$ equations of the system \eqref{L9}.

$(D)$ There exists a non-decreasing function $\sigma(z)$ having exactly $k$ nonzero points of growing, including $a$ but not including $b$, where $2k\le n$, i.e., the following  representation hold
\begin{equation}\label{L12}
c_j=\sum_{s=2}^k z_s^{j-1}\sigma_s+a^{j-1}\sigma_{k+1}, \ \ j=1,\ldots,n,
\end{equation}
where  \eqref{L3}, \eqref{L4} and 
\begin{equation}\label{L13} 
\sigma_2,\ldots,\sigma_{k+1}>0
\end{equation}
hold if and only if
\begin{enumerate}
\item[$(D_1)$] the matrices $\{c^{a}_{i+j-1+d}\}_{i,j=1}^{k}$ are singular for $0\le d\le n-2k$, i.e., 
\begin{equation}\label{L14} 
\det\{c^{a}_{i+j-1+d}\}_{i,j=1}^{k}=0, \ \ 0\le d\le n-2k;
\end{equation}

\item[$(D_2)$] the matrices $\{c^{a}_{i+j-1}\}_{i,j=1}^{k-1}$ and $\{c^{a}_{i+j+1}\}_{i,j=1}^{k-1}$ are positive definite;

\item[$(D_3)$] the matrix $\{c_{i+j-1}^b\}_{i,j=1}^{k}$ is positive definite.
\end{enumerate} 

The points $z_2,\ldots,z_k$ can be found as the roots of the equation
\begin{equation}\label{L5_4}
\det\left(\begin{array}{cccc}
c_1^a&c_2^a&\cdots&c_{k}^a\\
\cdots&\cdots&\cdots&\cdots\\
c^{a}_{k-1}&c^{a}_{k}&\cdots&c^{a}_{2k-2}\\
1&z&\cdots&z^{k-1}
\end{array}
\right)=0.
\end{equation}
If $z_2,\ldots,z_k$ are known, the numbers $\sigma_2,\ldots,\sigma_{k+1}$ can be found from the first $k$ equations of the system \eqref{L12}.
\end{lemma}

The proof of Lemma~\ref{Lem} is given in Appendix~\ref{app2}.

\section{Algorithm of finding the optimal control}\label{sec_3}

The analysis given in the previous section suggests the following algorithm of finding the optimal control: for a given point $x^0$, consider one by one all possible cases from Table~\ref{tab1} (four or five, in dependence on the parity of $n$) and check the appropriate part of Lemma~\ref{Lem}. If the conditions of Lemma~\ref{Lem} are satisfied, solve the corresponding equations; comparing all such cases, find the optimal time and optimal control. 

If none of these cases is appropriate, then the point $x^0$ does not belong to the 0-controllability domain, i.e., there does not exist a control satisfying the constraint $|u(t)|\le1$ and transferring the system \eqref{s1} from the point $x^0$ to the origin in a finite time. In fact, if such a control exists, then an optimal control exists as well due to the Filippov Theorem 
\cite{Fil1}, \cite{Fil2}. 
However, any optimal control relates to one of the cases described in Table~\ref{tab1}. On the other hand, for $n\ge3$ any neighborhood of the origin contains points that do not belong to the 0-controllability domain due to equations $\dot x_j=x_1^{j-1}$ for odd $j\ge3$. Thus, our algorithm allows us to determine if a given point $x^0$ belongs to the 0-controllability domain.

Let us analyze Table~\ref{tab1} more specifically. In each case, we define $c_1,\ldots,c_n$ by \eqref{s7}, \eqref{s9} and substitute $a=0$ and/or $b=x_1^0$ if appropriate. 

Trying Case 1, we substitute $a=0$ and $b=x_1^0$ in \eqref{s7}, \eqref{s9}, what gives
$$\begin{array}{l}
c_1=\widehat \theta-x_1^0,\\
\displaystyle c_{j}=-x_{j}^0-\frac{(x_1^0)^{j}}{j}, \ \ j=2,\ldots,n.
\end{array}
$$
Therefore, $c_1$ is a function on $\widehat\theta$ and $c_2,\ldots,c_n$ are known numbers. The conditions are of type $A$, i.e., they are described in part $(A)$ of Lemma~\ref{Lem} with  $k=m+1$. In this case we have the unique equation \eqref{L5} with $d=0$, therefore, we directly find $\widehat\theta$. Then we substitute $\widehat\theta$ to $c_1$ and check conditions $(A_2)$ and $(A_3)$. If they are satisfied, the solution exists. We find $z_2,\ldots,z_{k}$ as the roots of the equation \eqref{L5_1} and, finally, the numbers $\sigma_2,\ldots,\sigma_{k}$ are found from the system \eqref{L1}. 

In Cases 2 and 3, we substitute $a=0$, what gives
$$\begin{array}{l}
c_1=\widehat \theta+x_1^0-2b,\\
\displaystyle c_{j}=-x_{j}^0+\frac{(x_1^0)^{j}-2b^j}{j}, \ \ j=2,\ldots,n,
\end{array}
$$
i.e., $c_2,\ldots,c_n$ become polynomials in $b$ and $c_1$ is a polynomial in $b$ and $\widehat\theta$. In these cases we deal with types $A$ and $C$ respectively, where conditions \eqref{L5} and \eqref{L11} include two equations, with $d=0$ and $d=1$. However, the equation with $d=1$ does not include $c_1$ and, therefore, does not include $\widehat\theta$. Thus, this is a polynomial equation on $b$.  

Therefore, we find roots of the equation \eqref{L5} or \eqref{L11} with $d=1$ such that $0<b<x_1^0$ and, substituting each such a root to the equation with $d=0$, find  $\widehat\theta$. Then we substitute each obtained pair $(b,\widehat\theta)$ to $c_1,\ldots,c_n$ and check conditions $(A_2)$, $(A_3)$ (in Case~2) or  $(C_2)$, $(C_3)$ (in Case~3) to find out if the corresponding case is realized. If this is so, we find $z_2,\ldots,z_{k}$ using the equation \eqref{L5_1} or \eqref{L5_3}. For $\sigma_1,\ldots,\sigma_{k}$, we solve the system  \eqref{L1} or \eqref{L9}. 

Cases 4 and 7 are treated completely analogously, but, instead of $a$, we substitute $b=x_1^0$ in \eqref{s7}, \eqref{s9}, what gives
$$\begin{array}{l}
c_1=\widehat \theta-x_1^0+2a,\\
\displaystyle c_{j}=-x_{j}^0+\frac{-(x_1^0)^{j}+2a^j}{j}, \ \ j=2,\ldots,n,
\end{array}
$$
i.e., $c_2,\ldots,c_n$ become functions on $a$ and $c_1$  becomes a function on $a$ and $\widehat\theta$. Now we deal with types $A$ and $D$ respectively and also have two equations, for $d=0$ and $d=1$, in conditions  \eqref{L5} and \eqref{L14}. Analogously to Cases 2 and 3, first we find $a<0$ from the equation for $d=1$ and then we find $\widehat\theta$ from the equation for $d=0$. Afterwards, we check conditions $(A_2)$, $(A_3)$ (in Case~4) or  $(D_2)$, $(D_3)$ (in Case~7); if the corresponding case is realized, we find the solution using the equation \eqref{L5_1} or \eqref{L5_4} and the system  \eqref{L1} or \eqref{L12}. 

Let us pass to Cases 5, 6, 8, 9. Now, both variables $a$ and $b$ are unknown, i.e., $c_1,\ldots,c_n$ are defined by \eqref{s7}, \eqref{s9} without simplification. Then, $c_2,\ldots,c_n$ are polynomials in two variables $a$ and $b$, while $c_1$ is a polynomial in three variables $a$, $b$, and $\widehat\theta$. The corresponding conditions  \eqref{L5}, \eqref{L8}, \eqref{L11}, \eqref{L14} contain three equations, for $d=0$, $d=1$, and $d=2$. However, two latter equations (for $d=1$ and $d=2$) do not include  $\widehat\theta$. Therefore, we solve a system of two polynomial equations in two variables $a$ and $b$, and, substituting each a root to the equation for $d=0$, find  $\widehat\theta$. Then we substitute all obtained triples $a$, $b$, $\widehat\theta$ to check  conditions  $(A_2)$, $(A_3)$ (in Case~5),  $(C_2)$, $(C_3)$ (in Case 6),  $(D_2)$, $(D_3)$ (in Case~8),  $(B_2)$, $(B_3)$ (in Case 9) in order to find out if the case is realized. If this is so, we substitute $a$, $b$, $\widehat\theta$ to find the roots of the polynomial \eqref{L5_1}, \eqref{L5_2}, \eqref{L5_3}, or \eqref{L5_4} and solve the corresponding linear system \eqref{L1}, \eqref{L6}, \eqref{L9},  or \eqref{L12}.

We emphasize that, \emph{for any $n$, in the worst case a system of two polynomial equations in two variables should be solved}. For example, one can do that using a resultant of these two polynomials. 

Summarizing, we formulate our main result. 

\begin{theorem}\label{Thm}
The time-optimal control problem \eqref{s1}, \eqref{s2}, where $n\ge4$, is reduced to solving a truncated Hausdorff moment problem. For generic points with $x_1^0>0$, four (if $n$ is even) or five (if $n$ is odd) cases of the optimal control are possible. These cases are listed in Table~\ref{tab1}: cases 1, 3, 5, 7,  9 correspond to odd $n$ and cases 2, 4, 6, 8 correspond to even~$n$. 

In Case 1, the solution includes solving one linear equation for $\widehat\theta$, checking three matrices for positive definiteness, solving one polynomial equation in one variable for $z_i$,  and solving one linear system for $\sigma_i$.

In Cases 2, 3, 4, 7, the solution includes solving one polynomial equation for $a$ or $b$, one linear equation for $\widehat\theta$, checking three matrices for positive definiteness, solving one polynomial  equation in one variable for $z_i$,  and solving one  linear system  for $\sigma_i$.

In Cases 5, 6, 8, 9, the solution includes solving a system of two polynomial equations for $a$ and $b$, one linear equation for $\widehat\theta$, checking three matrices for positive definiteness, solving one polynomial  equation in one variable for $z_i$,  and solving one  linear system  for $\sigma_i$.
\end{theorem}

\section{Examples}\label{sec_4}

We illustrate our algorithm by the time-optimal control problem for the four-dimensional system 
\begin{equation}\label{E2}
\dot x_1=u, \ \ \dot x_2=x_1, \ \ \dot x_3=x_1^2, \ \ \dot x_4=x_1^3.
\end{equation}
Here $n=4$ and $m=2$. We suppose that our initial points are generic; at least, that it true when we deal with float point calculations.

\subsection{Initial point $x^0=(1,-2,-6,2)$.} We have four cases; let us first try Case 4. Here $k=2$ and the conditions are of type $A$ (see Tables~\ref{tab1} and~\ref{tab2}).

\begin{figure}
\centering
\includegraphics[width=3.5in]{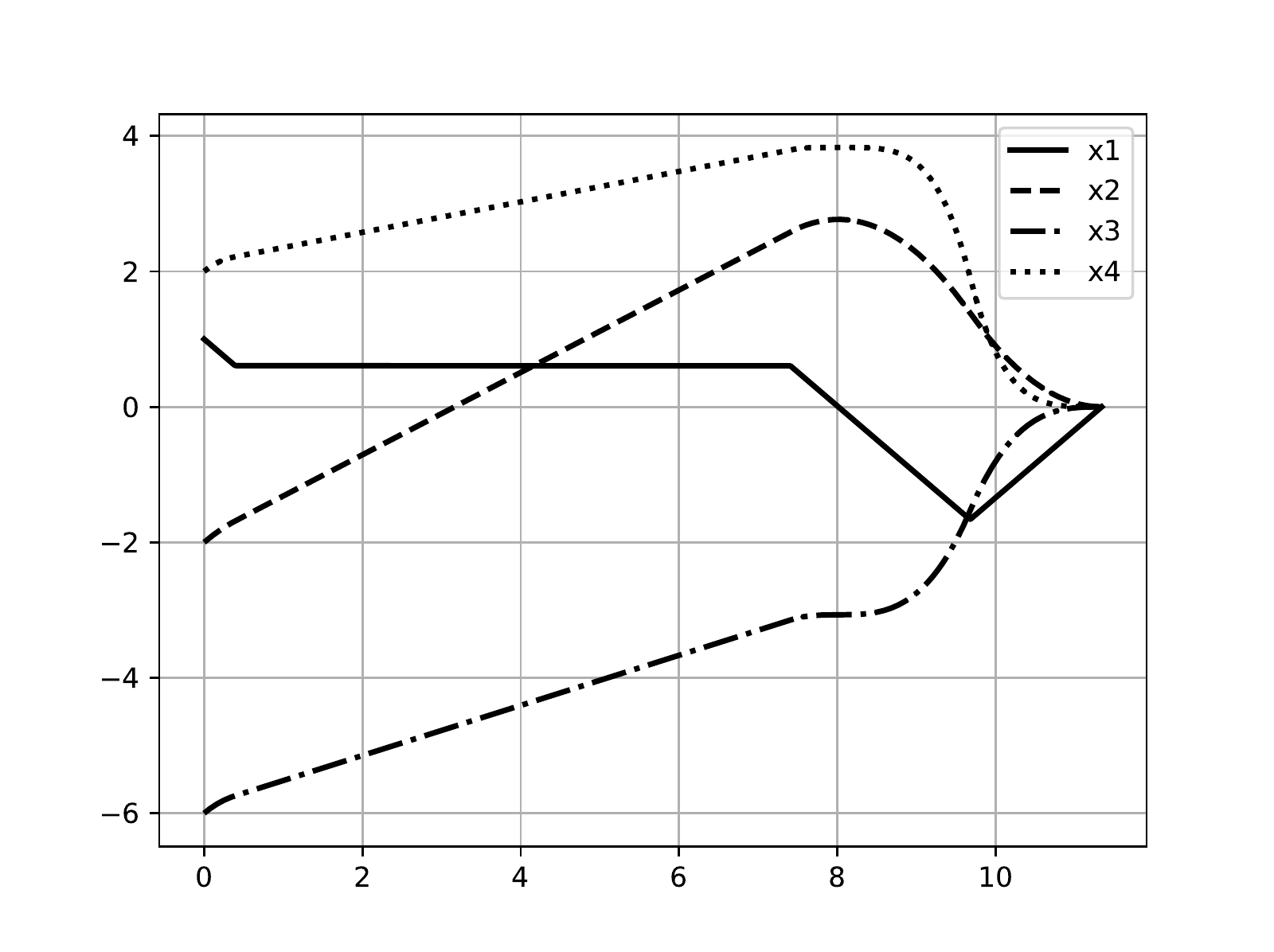}
\caption{Components of the optimal trajectory, $x^0=(1,-2,-6,2)$.}\label{fig1}
\end{figure}

First we use condition $(A_1)$ and solve equation \eqref{L5} with $d=1$ for $a$. Substituting $b=x_1^0=1$, we have
$$c_2=a^2+\frac32, \ \ c_3=\frac23a^3+\frac{17}3, \ \ c_4=\frac12a^4-\frac{9}4,
$$ 
therefore,
$$\begin{vmatrix}
c_2&c_3\\
c_3&c_4
\end{vmatrix}=\frac{1}{18}a^6 + \frac34a^4 - \frac{68}{9}a^3 - \frac94a^2 - \frac{2555}{72}.
$$
This polynomial has only one negative root $a\approx -1.66366$. Substituting it to $c_1$, $c_2$, and $c_3$, we get $c_1\approx \widehat\theta - 4.32731$, $c_2\approx 4.26776$, $c_3\approx2.59693$ and then find $\widehat\theta$ as a root of equation \eqref{L5} with $d=0$, 
$$\begin{vmatrix}
c_1&c_2\\
c_2&c_3
\end{vmatrix}=0,
$$
what gives $\widehat\theta\approx11.34087$. Now we substitute $\widehat\theta$ and check conditions $(A_2)$ and $(A_3)$: we have $c_1\approx7.01356>0$, $c_3\approx2.59693>0$, and $c_1^{a,b}=-c_3+(a+1)c_2-ac_1\approx6.23889>0$, therefore, they are satisfied. From \eqref{L5_1} we find $z_2$ as a root of the equation $c_1z-c_2=0$, what gives $z_2\approx 0.608501$. Finally, from the equation $c_1=\sigma_2$ we get $\sigma_2\approx 7.01356$. Components of the corresponding trajectory are shown in Fig.~\ref{fig1}.

Checking other cases, we see that there are no other possible candidates for optimal control. Therefore, the obtained control is optimal.

\subsection{Initial point $x^0=(1,-2,-\frac{12}{5},2)$.}  
\begin{figure}
\centering
\includegraphics[width=3.5in]{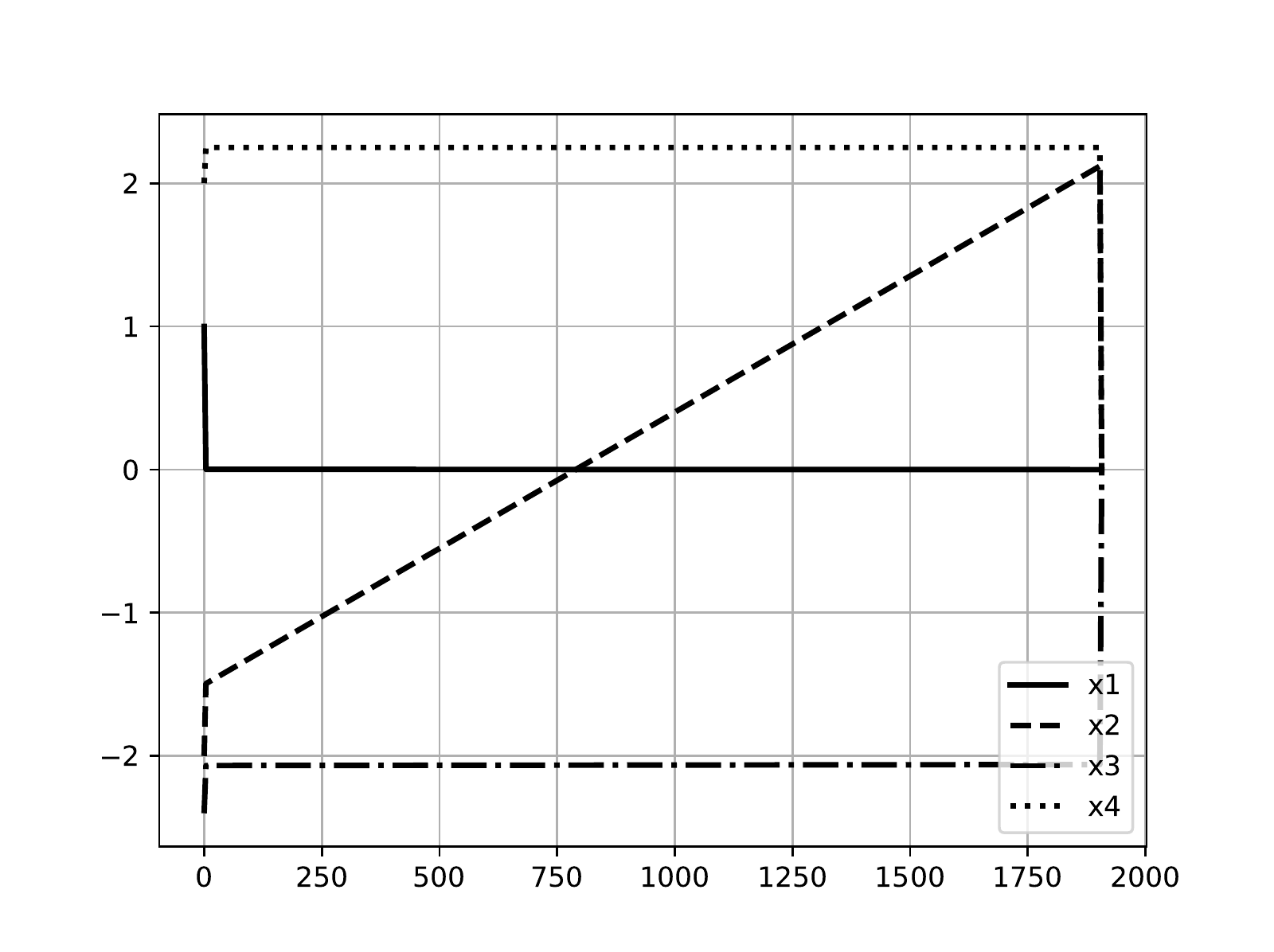}
\caption{Components of the optimal trajectory, $x^0=(1,2,-\frac{12}{5},\frac12)$.}\label{fig2}
\end{figure}

For this initial point, like the previous example, Case 4 is realized. However, now the optimal time is very large, $\widehat\theta\approx1907.10809$. Components of the optimal trajectory are shown in Fig.~\ref{fig2}. Here we have $z_2\approx 0.001903$ and $\sigma_2\approx1903.19513$. Thus, almost all the time the trajectory moves with very small first component $x_1(t)$, i.e., in a neighborhood of the ``equilibrium subspace'' $x_1=0$. This shows that $x^0$ is close to the boundary of the 0-controllability domain; one can prove that the point $(1,-2,-\frac{12}{5}+\frac{1}{100},2)$ cannot be steered to the origin by a control satisfying the constraint $|u|\le1$. 

\subsection{Initial point $x^0=(1,2,-3,\frac12)$.} As above, first try Case 4. Applying condition $(A_1)$ and substituting $b=x_1^0=1$, we obtain the polynomial
$$\begin{vmatrix}
c_2&c_3\\
c_3&c_4
\end{vmatrix}=\frac{1}{18}a^6 - \frac54a^4 - \frac{32}{9}a^3 - \frac34a^2 - \frac{377}{72};
$$
again it has only one negative root $a\approx -1.58127$. From the equation \eqref{L5} with $d=0$ we find $\widehat\theta\approx4.16254$; substituting $a$ and $\widehat\theta$ we get $c_1\approx 0.000005>0$, $c_3\approx0.030804>0$, hence, condition $(A_2)$ holds. However, $c_1^{a,b}=-c_3+(a+1)c_2-ac_1\approx-0.03103<0$, therefore, condition $(A_3)$ is not satisfied. Thus, this case is wrong. 

It can be shown that for the initial point under consideration Case 8 is suitable, then $k=1$ and conditions are of type $D$. Let us treat this case in detail.

First, condition $(D_1)$ means that we should solve a system of two polynomial equations with respect to two variables $a$ and $b$, namely, equations \eqref{L14} with $d=1$ and $d=2$. They are of the form
\begin{equation}\label{E1}
\begin{array}{l}
c_2^a=-\frac13a^3-\frac23b^3+ab^2+\frac32a+\frac{10}3=0,\\[5pt]
c_3^a=-\frac16a^4-\frac12b^4+\frac23ab^3-\frac{10}{3}a-\frac14=0.
\end{array}
\end{equation}
We may use the resultant of the polynomials $c_2^a$ and $c_3^a$. To this end, let us write them as polynomials of $b$ with coefficients being polynomials of $a$, 
$$\begin{array}{c}
c_2^a=-\frac23b^3+ab^2+(-\frac13a^3+\frac32a+\frac{10}3),\\[5pt]
c_3^a=-\frac12b^4+\frac23ab^3+(-\frac16a^4-\frac{10}{3}a-\frac14);
\end{array}$$
we seek for such values of $a$ for which these polynomials have a common root. As is well known, such values are roots of the determinant of the Sylvester matrix,
$$\begin{vmatrix}
-\frac23&a&0&e_1&0&0&0\\
0&-\frac23&a&0&e_1&0&0\\
0&0&-\frac23&a&0&e_1&0\\
0&0&0&-\frac23&a&0&e_1\\
-\frac12&\frac23a&0&0&e_2&0&0\\
0&-\frac12&\frac23a&0&0&e_2&0\\
0&0&-\frac12&\frac23a&0&0&e_2\\
\end{vmatrix}=0,
$$
where we use the notation $e_1=-\frac13a^3+\frac32a^2+\frac{10}3$, $e_2=-\frac16a^4-\frac{10}3a-\frac14$. This polynomial equation (of degree 8) has four negative solutions; we substitute each of them to $c_2^a$ and $c_3^a$, obtain polynomials of $b$, and find their common roots. This gives all solutions $a,b$ of the system \eqref{E1}. In our case we have only two pairs $a,b$ such that $a<0$ and $b>x_1^0=1$. For the first pair, $a\approx-2.26126$ and $b\approx1.12296$, we obtain $\widehat\theta\approx4.728202$, however, $c_1^b\approx-3.52041<0$, i.e., condition $(D_3)$ does not hold. For the second pair, $a\approx-0.71829$ and $b\approx1.240801$, we obtain $\widehat\theta\approx6.43157$ and   $c_1^b\approx6.88305>0$, hence, condition $(D_3)$ is satisfied (since $k=1$, condition $(D_2)$ disappears). Finally, we find $\sigma_1=c_1\approx3.51338$.

\begin{figure}
\centering
\includegraphics[width=3.5in]{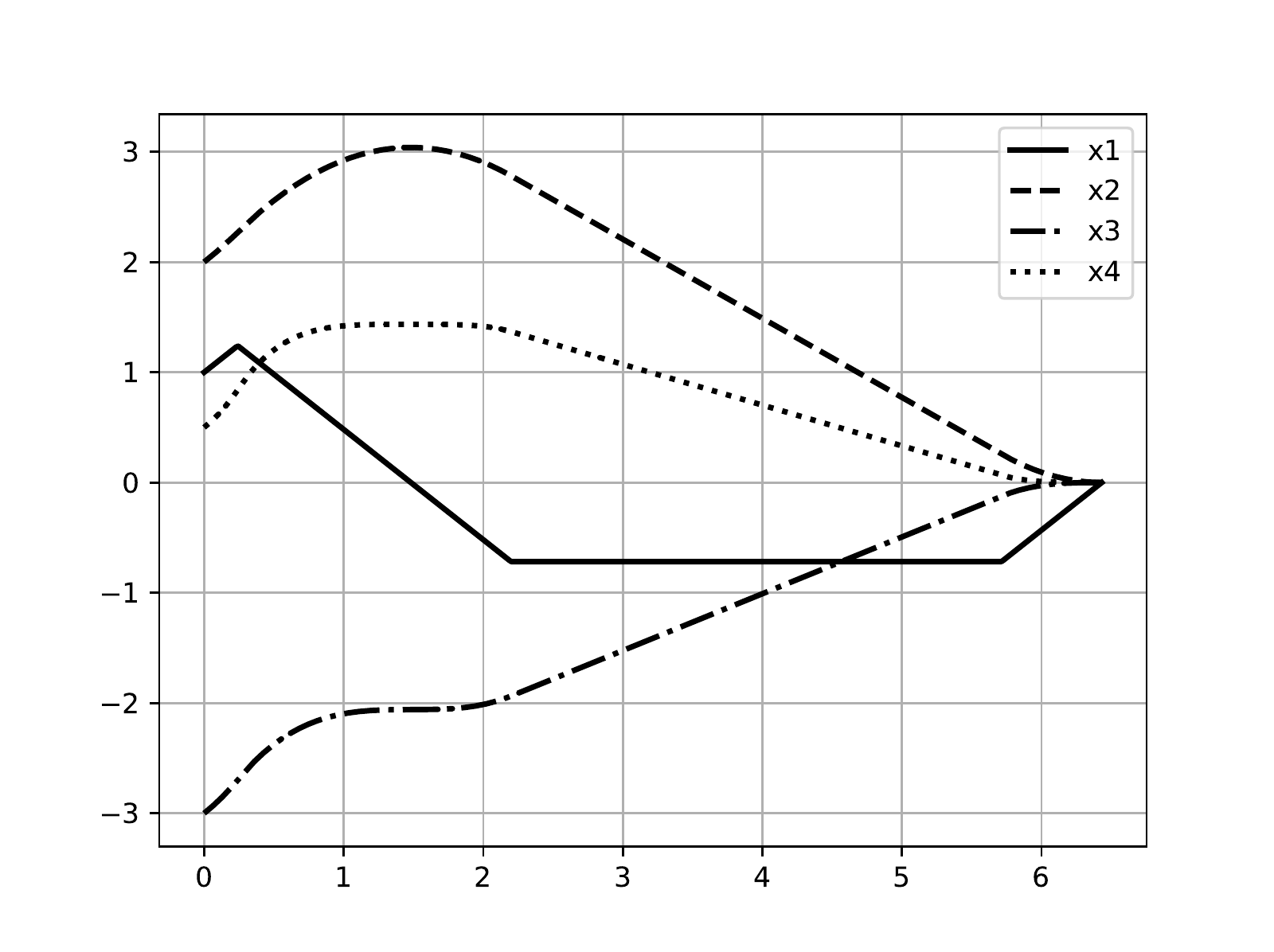}
\caption{Components of the optimal trajectory, $x^0=(1,2,-3,\frac12)$.}\label{fig3}
\end{figure}

Since there are no other possible candidates, we conclude that the obtained control is optimal. Components of the optimal trajectory are shown in Fig.~\ref{fig3}.

\subsection{Initial points $x^0=(1, -8, -3.8289, -1.8792)$ and $x^0=(1,-8,-28.4649,-1.8792)$.} 
\begin{figure}
\centering
\includegraphics[width=3.5in]{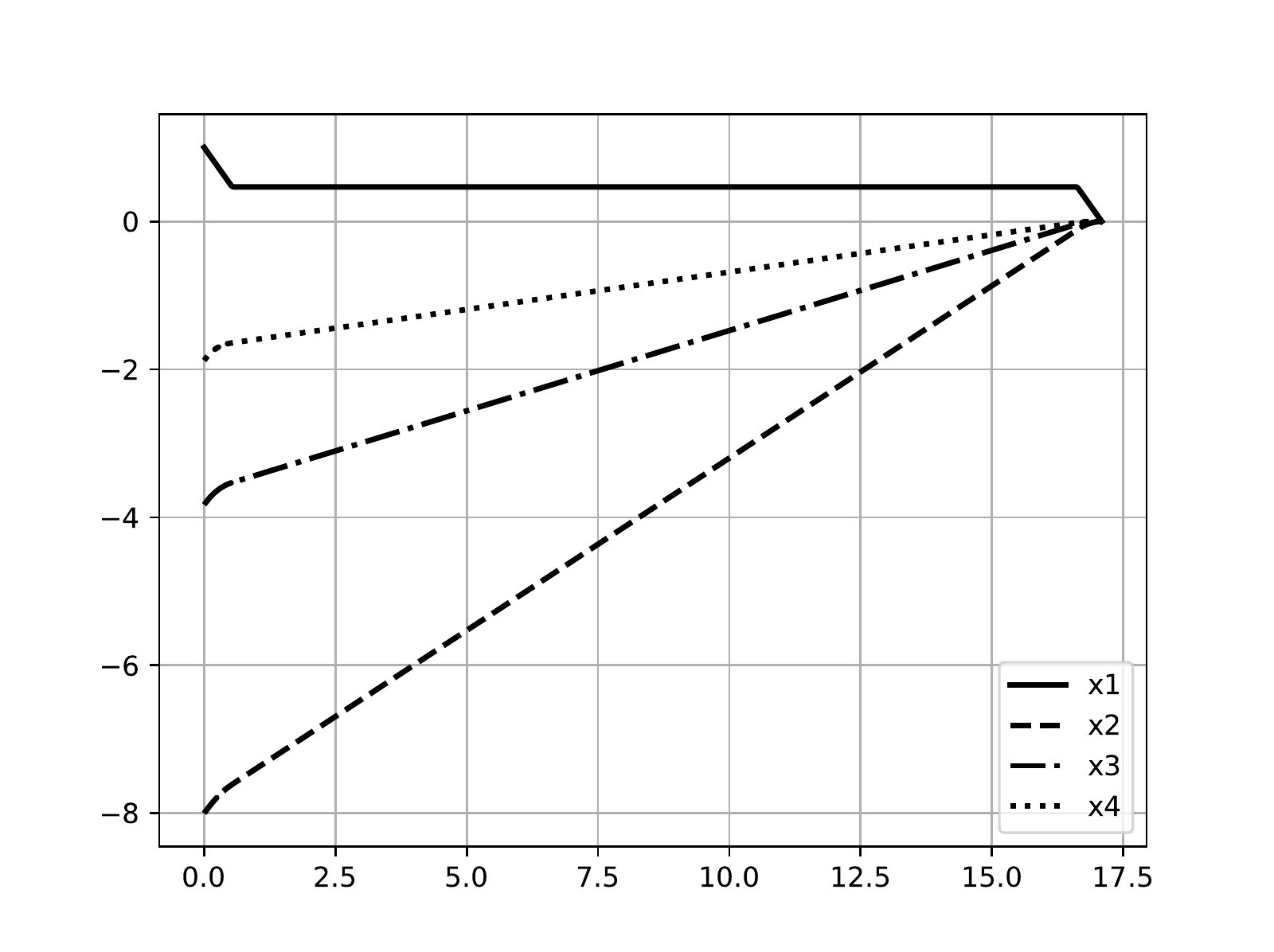}
\caption{Components of the optimal trajectory, ${x^0=(1, -8, -3.8289, -1.8792)}$.}\label{fig4}
\end{figure}

\begin{figure}
\centering
\includegraphics[width=3.5in]{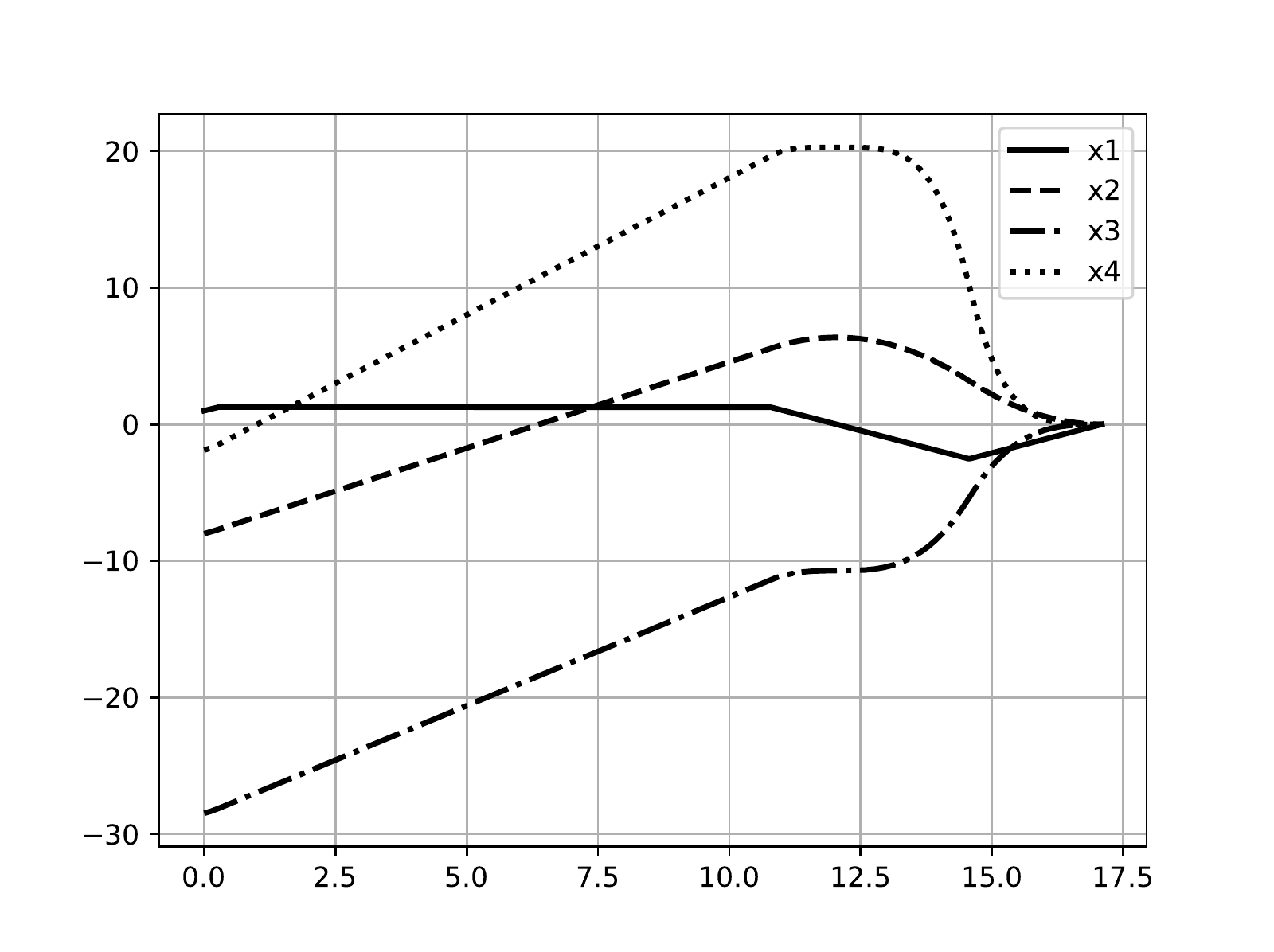}
\caption{Components of the optimal trajectory, ${x^0=(1,-8,-28.4649,-1.8792)}$.}\label{fig5}
\end{figure}

In 
\cite{I} 
the three-dimensional  time-optimal control problem for the system 
\begin{equation}\label{E3}
\dot x_1=u, \ \ \dot x_2=x_1, \ \ \dot x_3=x_1^3
\end{equation}
was completely solved. Such a problem can be thought of as a time-optimal control problem for a four-dimensional system \eqref{E2} where the coordinate $x^0_3$ is free. This feature leads to more complicated structure of optimal controls (though, equations for finding optimal controls are easier); in particular, an optimal control can be one of eight possible types. Moreover, some initial points admit two optimal controls of different types. For example, let us consider the line in the four-dimensional space such that $x_1=1$, $x_2=-8$, $x_4\approx -1.8792$ (the precise value can be found in 
\cite{I})
and $x_3$ can be arbitrary. One can show that the optimal time $\widehat\theta$ as a function of $x_3$ has \emph{two minimum points with the same minimum value}, namely, $x_3\approx-3.8289$ and $x_3\approx-28.4649$; the minimum value is $\widehat\theta\approx17.0918$. The point $x^0=(1, -8, -3.8289, -1.8792)$ is not generic in the sense mentioned above: the optimal control corresponds to Case 2 with $b=x_1^0$ or, what is the same, to Case 4 with $a=0$. The optimal control for the point   $x^0=(1, -8, -28.4649, -1.8792)$ is of Case~6. So, we obtain two different solutions for the system \eqref{E3} found in 
\cite{I};
components of the optimal trajectories are shown in Fig.~\ref{fig4} and Fig.~\ref{fig5}. The obtained result shows, in particular, that controllability sets of the system \eqref{E2} can be non-convex.

\section{Conclusion and open questions}

In the paper, the time-optimal problem for dual-to-integrator systems \eqref{s1}--\eqref{s2} was considered. We showed that this problem is reduced to a truncated Hausdorff moment problem. We proved that, regardless of system's dimension, only four or five cases (in dependence on dimension's parity) of the Hausdorff moment problem arise. In these moment problems, the interval $[a,b]$ where the problem is stated can be unknown. In each case one needs to solve a system of at most two polynomial equations in two variables $a$ and $b$ and, after that, a polynomial equation in one variable to find a control (in some cases, one solves only one or two polynomial equations in one variable). As a result, we obtain an analytic solution of the time-optimal problem \eqref{s1}--\eqref{s2}, which is in principle the same for all dimensions.

Finally, we formulate \emph{open questions and perspectives} concerning the time-optimal control problem \eqref{s1}--\eqref{s2}. 
\begin{enumerate}
\item As was mentioned in Section~\ref{sec_3}, the 0-controllability domain of the system \eqref{s1} does not coincide with the whole space ${\mathbb{R}}^n$ for $n\ge3$. In this context, realization of one of the cases from Table~\ref{tab1} means that a given point $x^0$ belongs to the 0-controllability domain. However, \emph{the problem of an explicit analytic description} of the 0-controllability domain remains unsolved except of the simplest case $n=3$ studied in 
\cite{SIAM2003}. 
The general case $n\ge4$ is a topic for further research.
\item 
Each case from Table~\ref{tab1} defines a set of points $x^0$ for which this case is realized, i.e., the corresponding control exists. A union of these sets (and their boundaries corresponding to non-generic points) is a 0-controllability domain. A question arises \emph{if these sets can overlap}. In the general case it is possible, what was shown in 
\cite{I}.
On the other hand, examples considered above in Section~\ref{sec_4} suggest that for systems \eqref{s1} this is not true. 
\item If the sets from the previous question \emph{can} overlap, the \emph{uniqueness of the optimal control} should be studied. 
\item Within several cases of our algorithm, when finding the optimal time, one solves a system of two polynomial equations in two variables. A resultant of these polynomials can be used, however, its degree is very high. The question is if another polynomial in one variable can be constructed having a degree less than a resultant. Alternatively, a method of effective numerical solving could be proposed, which uses a particular form of the polynomials under consideration.
\item As was mentioned in Section~\ref{sec1}, the system \eqref{s1} approximates a certain class of affine systems in the sense of time optimality, namely, systems that after some change of variables take the form
\begin{equation}\label{eq_e}
\begin{array}{l}
\dot x_1=u+p_1(x)+q_1(x)u,\\
\dot x_j=x_1^{j-1}+p_j(x)+q_j(x)u, \ j=2,\ldots,n,
\end{array}
\end{equation}
where $p_j(x)$, $q_j(x)$ are real analytic functions whose Taylor series contain terms $x_1^{k_1}\cdots x_n^{k_n}$ such that $k_1+2k_2+\cdots+nk_n\ge j$, $j=1,\ldots,n$. This means that the optimal times and optimal controls for systems  \eqref{s1} and \eqref{eq_e} are equivalent as $x^0\to0$. A perspective direction is to develop a method of successive approximation for solving the time-optimal problem for systems \eqref{eq_e}, which deals with the time-optimal problem for the system \eqref{s1} on each step, like the case of linear systems studied in 
\cite{JMAA}. 

\end{enumerate}
\appendix
\begin{appendices}

\section{Truncated Hausdorff moment problem}\label{app1}

In this appendix we recall some basic results on a truncated Hausdorff moment problem, which are used in the present paper; an extensive exposition and discussions can be found in the book
\cite{KN}.

The \emph{truncated Hausdorff moment problem on the interval $[a,b]$} can be formulated as follows: given a sequence of numbers $c_1,\ldots,c_n$, 

--- find out if there exists a non-decreasing function $\sigma(z)$ such that the following \emph{moment equalities} hold
\begin{equation}\label{s10_1-1}
c_j=\int_a^bz^{j-1}d\sigma(z), \ \ j=1,\ldots, n;
\end{equation}

--- if this is the case, determine if a solution is unique;

--- if so, find the solution; if not, find solutions of the simplest form.

Such a problem is also referred to as a \emph{mass distribution problem} due to its mechanical interpretation: a function $\sigma(z)$ defines a mass distribution of a rod whose ends have coordinates $a$ and $b$, namely, $\sigma(z)$  equals the mass of the piece between points $a$ and $z$. Then the problem consists in finding a mass distribution whose $n$ first moments are given; in particular, $c_1=\int_a^bd\sigma(z)$ equals the mass of the whole rod and $c_2=\int_a^bzd\sigma(z)$, $c_3=\int_a^bz^2d\sigma(z)$ define the coordinate of its center of mass and the moment of inertia about an axis at the point $z=0$.

Necessary and sufficient conditions of solvability of the Hausdorff moment problem \eqref{s10_1-1} are based on the Markov-Luk\'{a}cs Theorem on non-negative polynomials 
\cite{L}, \cite{KN}.
Namely, the problem \eqref{s10_1-1} is solvable iff the following two matrices are non-negative definite:

--- for $n=2m+1$ 
\begin{equation}\label{s10_1-2}
\{c_{i+j-1}\}_{i,j=1}^{m+1} \  \mbox{ and } \  \{-c_{i+j+1}+(a+b)c_{i+j}-abc_{i+j-1}\}_{i,j=1}^{m};
\end{equation}

--- for $n=2m$ 
\begin{equation}\label{s10_1-3}
\{c_{i+j}-ac_{i+j-1}\}_{i,j=1}^{m}  \ \mbox{ and } \  \{-c_{i+j}+bc_{i+j-1}\}_{i,j=1}^{m}.
\end{equation}

Among possible solutions, the class of step functions $\sigma(z)$ is considered separately. A step function describes masses that are concentrated at its points of discontinuity. Namely, if $b\ge z_1>\ldots>z_{k+1}\ge a$ are points of discontinuity of a step function $\sigma(z)$  and $\sigma_1,\ldots,\sigma_{k+1}>0$ are jump values at these points, then moment equalities \eqref{s10_1-1} take the form
\begin{equation}\label{s10_1-4}
c_j=\sum_{s=1}^{k+1} z_s^{j-1}\sigma_s, \ \ j=1,\ldots,n.
\end{equation}
Each such function can be assigned with an \emph{index} defined as follows: any point of the open interval $(a,b)$ is assumed to have index 2, while points $a$ and $b$ are assumed to have index~1; the index of a step function is defined as the sum of indexes of all its points of discontinuity. The mentioned solvability theorem can be supplemented as follows. If both matrices in \eqref{s10_1-2} (for odd $n$) or in \eqref{s10_1-3} (for even $n$) are positive definite, then there exists an infinite set of solutions, among which there exist exactly two solutions being step functions of index $n$. If both matrices are non-negative and at least one of them is singular, the solution is unique; it is a step function of index no greater than $n-1$.

Geometrically, the solvability set for the problem \eqref{s10_1-1} is a conic hull of the curve $\{(1,t,\ldots,t^{n-1}): a\le t\le b\}$. The case of positive definiteness of the matrices \eqref{s10_1-2} or \eqref{s10_1-3} means that the point $c=(c_1,\ldots,c_n)$ belongs to the interior of the solvability set, while if one of the matrices is singular, then the point lies on the boundary of the solvability set. In the latter case the points of discontinuity of $\sigma(z)$ can be found as roots of a polynomial explicitly expressed via $c_1,\ldots,c_n$. For example, if $n=2m+1$ and the matrix $\{c_{i+j-1}\}_{i,j=1}^{m+1}$ is singular, the points of discontinuity of $\sigma(z)$ are the roots of the polynomial
$$
\det\left(\begin{array}{cccc}
c_1&c_2&\cdots&c_{m+1}\\
\cdots&\cdots&\cdots&\cdots\\
c_{m}&c_{m+1}&\cdots&c_{2m}\\
1&z&\cdots&z^{m}
\end{array}
\right)=0.
$$
When the points of discontinuity are known, the jump values are found from the system of \emph{linear} equations \eqref{s10_1-4}.

If the index of the solution is less than $n-1$, then matrices of smaller dimension are singular. Conditions given in Lemma~\ref{Lem} describe such a case with the additional requirement that $z=0$ is not a point of discontinuity of $\sigma(z)$. 
 
\section{Proof  of Lemma~\ref{Lem}.}\label{app2}

We prove case $(A)$ only; the other cases can be proved completely analogously.  

\emph{Necessity}. Suppose \eqref{L1}--\eqref{L2} hold. Denote by $\widehat Q(z)$ a nontrivial polynomial with the roots $z_2,\ldots,z_k$,  
$$\widehat Q(z)=\prod_{s=2}^k(z-z_s)=\sum_{i=1}^k\widehat q_iz^{i-1}.
$$ 
Substituting \eqref{L1}, for any $0\le d\le n+1-2k$ and for any $j=1,\ldots,k$ one obtains 
$$\sum_{i=1}^kc_{i+j-1+d}\widehat q_i=\sum_{i=1}^k\sum_{s=2}^kz_s^{i+j-2+d}\sigma_s\widehat q_i=
$$
$$=\sum_{s=2}^k\left(\sum_{i=1}^k\widehat q_iz_s^{i-1}\right)z_s^{j-1+d}\sigma_s=\sum_{s=2}^k\widehat Q(z_s)z_s^{j-1+d}\sigma_s=0.
$$
These equalities can be rewritten in the matrix form, 
$$\left(\begin{array}{cccc}
c_{1+d}&c_{2+d}&\cdots&c_{k+d}\\
c_{2+d}&c_{3+d}&\cdots&c_{k+1+d}\\
\ldots&\ldots&\ldots&\ldots\\
c_{k+d}&c_{k+1+d}&\cdots&c_{2k-1+d}
\end{array}\right)\left(\begin{array}{c}
\widehat q_1\\\widehat q_2\\\ldots\\\widehat q_{k}
\end{array}\right)=0,
$$
where $0\le d\le n+1-2k$, what proves $(A_1)$. Moreover, the equality 
$$\left(\begin{array}{cccc}
c_{1}&c_{2}&\cdots&c_{k}\\
\ldots&\ldots&\ldots&\ldots\\
c_{k-1}&c_{k}&\cdots&c_{2k-2}\\
1&z&\ldots&z^{k-1}
\end{array}\right)\left(\begin{array}{c}
\widehat q_1\\\widehat q_2\\\ldots\\\widehat q_{k}
\end{array}\right)=0
$$
holds for $z=z_2$, \ldots, $z=z_k$, hence, $z_2,\ldots,z_k$ are the roots of the equation \eqref{L5_1}. 

To prove $(A_2)$, let us consider a quadratic form 
$$\sum_{i,j=1}^{k-1}c_{i+j-1}q_iq_j=\sum_{i,j=1}^{k-1}\sum_{s=2}^kz_s^{i+j-2}\sigma_sq_iq_j=
$$
\begin{equation}\label{Lp_0}
=\sum_{s=2}^{k}\left(\sum_{i=1}^{k-1}q_iz_s^{i-1}\right)^2\sigma_s=\sum_{s=2}^kQ^2(z_s)\sigma_s,
\end{equation}
where $Q(z)=\sum_{i=1}^{k-1}q_iz^{i-1}$ is of degree no more than $k-2$. If $Q(z)$ is nontrivial, then all terms in the right hand side of \eqref{Lp_0} are non-negative and no more than $k-2$ of them vanish (we take into account \eqref{L2}). Therefore, at least one term is positive, hence, the sum in the right hand side of \eqref{Lp_0} is positive. This implies that the matrix $\{c_{i+j-1}\}_{i,j=1}^{k-1}$ is positive definite. For the matrix $\{c_{i+j+1}\}_{i,j=1}^{k-1}$, the same considerations give 
$$\sum_{i,j=1}^{k-1}c_{i+j+1}q_iq_j=\sum_{s=2}^kQ^2(z_s)z_s^2\sigma_s,
$$ 
where $Q(z)$ is of degree no more than $k-2$. Arguing analogously and taking into account \eqref{L4}, we conclude that the matrix $\{c_{i+j+1}\}_{i,j=1}^{k-1}$ is positive definite, what implies $(A_2)$.

Analogously, 
$$\sum_{i,j=1}^{k-1}c^{a,b}_{i+j-1}q_iq_j=$$
$$=\sum_{i,j=1}^{k-1}(-c_{i+j+1}+(a+b)c_{i+j}-abc_{i+j-1})q_iq_j=
$$ 
$$=\sum_{s=2}^{k}(z_s-a)(b-z_s)Q^2(z_s)\sigma_s,
$$
where $Q(z)$ is of degree no more than $k-2$. Taking into account \eqref{L3}, for a nontrivial $Q(z)$ we get that this sum is positive, what proves $(A_3)$.

\emph{Sufficiency.} Suppose conditions  $(A_1)$--$(A_3)$ hold. Let us consider the Hausdorff moment problem for the sequence $\{c_1,\ldots,c_{2k-1}\}$, i.e.,
\begin{equation}\label{Lp_1}
c_j=\int_a^bz^{j-1}d\sigma(z), \ \ j=1,\ldots, 2k-1.
\end{equation}
Now we make use of the solvability conditions for this moment problem 
\cite[Chapter III]{KN}.
Due to suppositions  $(A_2)$ and $(A_1)$ for $d=0$, the matrices $\{c_{i+j-1}\}_{i,j=1}^{k-1}$ and $\{c_{i+j+1}\}_{i,j=1}^{k-1}$ are positive definite, while the matrix $\{c_{i+j-1}\}_{i,j=1}^{k}$ is singular. Hence, the matrix $\{c_{i+j-1}\}_{i,j=1}^{k}$ is non-negative definite. Taking into account also condition $(A_3)$, we obtain that the problem \eqref{Lp_1} has a unique solution $\sigma(z)$ which is a step function with exactly $k-1$ points of growing inside the interval $(a,b)$. Now, suppose $b>z_2>\cdots>z_{k}>a$ are points of growing of $\sigma(z)$ and $\sigma_2,\ldots,\sigma_{k}>0$ are the corresponding jump values, i.e., 
$$c_j=\sum_{s=2}^kz_s^{j-1}\sigma_s, \ \ j=1,\ldots,2k-1.
$$
Then, analogously to the arguments mentioned above, we get
\begin{equation}\label{Lp_3}
\sum_{i,j=1}^{k-1}c_{i+j+1}q_iq_j=\sum_{s=2}^kQ^2(z_s)z_s^2\sigma_s,
\end{equation}
where  $Q(z)=\sum_{i=1}^{k-1}q_iz^{i-1}$ is of degree no more than $k-2$. Suppose some jump point equals zero, $z_{s_0}=0$. Then there exists a nontrivial polynomial of degree $k-2$ with the roots $\{z_s: s\ne s_0\}$. Then the right hand side of \eqref{Lp_3} equals zero, what is impossible since $\{c_{i+j+1}\}_{i,j=1}^{k-1}$ is positive definite. Hence, all jump points $z_s$ are nonzero. 

Thus, the moment equalities \eqref{L1} hold for $j=1,\ldots,2k-1$ and \eqref{L3}--\eqref{L2} are satisfied. Let us prove that \eqref{L1} hold also for $q=2k,\ldots,n$ (if $2k\le n$). To this end, we use the rest of conditions $(A_1)$.

As one can prove 
\cite[Part VII, \S2, Problem 20]{PS},
conditions $(A_1)$ imply that $\det\{c_{i+j+d}\}_{i,j=1}^{k-1}$ for $d=0,\ldots,n+1-2k$ equal zero or do not equal zero simultaneously. By supposition $(A_2)$, $\det\{c_{i+j+1}\}_{i,j=1}^{k-1}\ne0$, therefore, $\det\{c_{i+j+d}\}_{i,j=1}^{k-1}\ne0$ for $d=0,\ldots,n-2k$. This means that $c_{2k},\ldots,c_n$ can be uniquely expressed via $c_1,\ldots,c_{2k-1}$ one by one from the equalities $\det\{c_{i+j-1+d}\}_{i,j=1}^k=0$ successively considered for $d=1,\ldots,n+1-2k$ (in fact, when expanding the determinant $\det\{c_{i+j-1+d}\}_{i,j=1}^k$ by the last row, $c_{2k-1+d}$ is multiplied by $\det\{c_{i+j-1+d}\}_{i,j=1}^{k-1}\not=0$). 

Now let us denote 
$$c'_j=\sum_{s=2}^k z_s^{j-1}\sigma_s, \ \ j=2k,\ldots,n.
$$ 
If we substitute $c_j$ by $c_j'$, $j=2k,\ldots,n$, in the determinants $\det\{c_{i+j-1+d}\}_{i,j=1}^k$, these determinants still vanish (as was shown above in the ``necessity'' part). Hence,  $c_{2k}',\ldots,c'_{n}$ are expressed via $c_1,\ldots,c_{2k-1}$ by the same formulas as  $c_{2k},\ldots,c_{n}$. Therefore, $c_j=c_j'$ for $j=2k,\ldots,n$, what proves \eqref{L1} for $j=2k,\ldots,n$. \qed

\end{appendices}

\section*{Acknowledgment} The authors would like to thank Sergey~Shugarev for numerous discussions on the matter of the present research.

\end{document}